\documentclass[12pt]{article}
\usepackage{graphicx}        

\usepackage{amssymb, amsmath}
\usepackage{bbm}

\usepackage{psfrag}
\usepackage{epsfig}

\usepackage{tikz}

\usepackage{latexsym}
\usepackage{amscd}
\newtheorem{theorem}{Theorem}

\def\B{{\cal B}}

\def\R{{\mathbb R}}
\def\Rd{\R^d}

\def\PROB{{\mathbb P}}
\def\EXP{{\mathbb E}}

\newcommand{\IND}{\mathbbm{1}}

\def\qed{\hfill {\ \vrule width 1.5mm height 1.5mm \smallskip}}

\def\be{\begin{equation}}
\def\ee{\end{equation}}

\def\qedskip{\smallskip\noindent}
\def\qed{\hfill $\Box$ \qedskip}              

\def\ol{\overline}
\def\defeq{\stackrel{\text{def}}{=}}

\def\diam{\text{diam}}

\begin{document}
\begin{titlepage}
\thispagestyle{empty}

\title{On the measure of Voronoi cells
\thanks{L\'aszl\'o Gy\"orfi was supported in part by the National Development Agency 
(NF\"U, Hungary) as part of the project Introduction of Cognitive
Methods for UAV Collision Avoidance Using Millimeter Wave Radar, grant
KMR-12-1-2012-0008.
Luc Devroye was supported by the Natural Sciences and Engineering
Research Council (NSERC) of Canada.
G\'abor Lugosi was supported by the Spanish Ministry of Science and Technology grant MTM2012-37195.
}
}
\author{Luc Devroye
\thanks{McGill University} \and 
 \and L\'aszl\'o Gy\"orfi
\thanks{Budapest University of Technology and Economics}
 \and 
G\'abor Lugosi\thanks{ICREA and Pompeu Fabra University}
\and Harro Walk\thanks{Universit\"at Stuttgart}}

\maketitle

\begin{abstract}
$n$ independent random points drawn from a density $f$ in $\Rd$ define a random
Voronoi partition. We study the measure of a typical cell of the partition. We prove
that the asymptotic distribution of the probability measure of the
cell centered at a
point $x \in \Rd$ is independent of $x$ and the density $f$. 
We determine all moments of
the asymptotic distribution and show that the distribution becomes more concentrated as $d$ becomes large. In particular, we show that the variance converges to zero
exponentially fast in $d$. 
We also obtain a density-free bound for the rate of convergence of the
diameter of a typical Voronoi cell.
\end{abstract}

\end{titlepage}

\section{Introduction}

Let $X_1,\ldots,X_n$ be independent, identically distributed random
vectors  taking values in $\R^d$. We denote the common distribution of
the $X_i$ by $\mu$. We assume throughout the paper that $\mu$ is
absolutely continuous with respect to the Lebesgue measure $\lambda$
and denote the density of $\mu$ by $f$. Hence, $\mu(A)=\int_A f(x)dx$
for all Lebesgue measurable sets $A\subset \Rd$.

The $X_i$ define a random partition of $\Rd$ into $n$ sets
$S_1,\ldots,S_n$ such that
$S_i$ contains all points in $\R^d$ whose nearest
neighbor among $X_1,\ldots,X_n$ is $X_i$. Ties are broken in favor
of smaller indices. (Because of the assumption of absolute continuity
of $\mu$, the tie-breaking rule is irrelevant throughout the paper.)
Formally,
\[
   S_i=\left\{x\in \Rd: \|x-X_i\|= \min_{j=1,\ldots,n} \|x-X_j\|
   \right\}  \bigcap \left\{x\in \Rd: \|x-X_i\| < \min_{j=1,\ldots,i-1} \|x-X_j\|
   \right\}~.
\]
$\{S_1,\ldots,S_n\}$ is a so-called \emph{Voronoi partition} and the
$S_i$ are the \emph{Voronoi cells}.

In this paper we are interested in the measure of a ``typical'' 
Voronoi cell. In particular, 
we study the conditional distribution of the random variable
$\mu(S_1)$ conditioned on the event that $X_1=x$ for some
$x$ in the support of $\mu$.


Note that since
\[
\sum_{j=1}^n \mu(S_j)=1
\]
and $\mu(S_1),\ldots ,\mu(S_n)$ are identically distributed, we have
\[
n\EXP\left[\mu(S_1))\right]= 1~.
\]
In Theorem \ref{thm:main} below we prove that, for $\mu$-almost
all $x$, we have $n\EXP\left[\mu(S_1)|X_1=x\right] \to 1$.
We also show that  $n^2\EXP\left[\mu(S_1)^2|X_1=x\right]$ converges
to a limit that is independent of $x$ and the distribution $\mu$. 
In fact, we prove that for $\mu$-almost
all $x$, $n\EXP\left[\mu(S_1)|X_1=x\right]$ has a limiting distribution
that only depends on the dimension. We show that the limiting
distribution becomes more concentrated as the dimension $d$ grows.

Finally, we study the diameter $\diam(S_1)$ of the Voronoi cell
centered at $X_1$. We show that for $\mu$-almost all $x$,
conditionally on $X_1=x$, $\diam(S_1)$ converges to zero at a rate
of $n^{-1/d}$.

Throughout the paper,
$B_{x,r}$ denotes the closed ball of radius $r>0$ centered at $x\in \R^d$.

\subsubsection*{Related work}

The measure of a ``typical'' cell in a Voronoi tessellation has been mostly
studied in the case when the points are drawn from a homogeneous Poisson
process. 
Asymptotically, this is equivalent to the special case 
of uniform distribution $\mu$ on (say) the unit ball. 
The study of the measure of Voronoi cells dates back to at least
Gilbert \cite{Gil62} who derived formulas and numerical estimates
for the second and third moments the measure of a Voronoi cell 
when $d=2$ or $3$. See also Brakke \cite{Bra85a}, \cite{Bra85b}, Hayen and Quine \cite{HaQu02}.

Our notion of the distribution of a typical cell is analogous
to the so-called ``Palm distribution'' of the volume of a Voronoi cell
in stochastic geometry---Stojan, Kendall, and Mecke \cite{StKeMe87}, M{\o}ller \cite{Mol94}, M{\o}ller and Stoyan \cite{MoSt07}.

Brakke \cite{Bra85a}, \cite{Bra85b},
Hayen and Quine \cite{HaQu02},
Heinrich et al. \cite{HeKoMeMu98}, 
Heinrich and Muche \cite{HeMu08},
Zuyev \cite{Zuy92}, and others
study characteristics of
``typical'' cells in a Voronoi tessallation of a homogeneous Poisson
process, including the second moment of the volume. 

For a survey and comprehensive treatment of
Voronoi diagrams, we refer to Okabe, Boots and Sugihara \cite{OkBoSu92} and
Okabe, Noots, Sugihara and Nok Chiu \cite{OkBoSuNo00}.

\section{Results}

Theorem \ref{thm:main} below establishes the asymptotic value
of the first and second moments of the measure of a typical cell centered at
a point $x$. The remarkable feature is that the asymptotic values
are independent of both the density $f$ and the point $x$ (for $\mu$-almost all $x$) and only
depend on the dimension $d$. In fact, in Theorem \ref{thm:Luc3} we
show that the limit distribution is also independent of $f$ and $x$.
We emphasize that both
theorems hold without any assumption on the density $f$.

The asymptotic second moment is expressed in terms of a random
variable $W$ defined as follows. 
Let $Y$ be a random vector uniformly distributed in $B_{0,1}$. Define
$\ol{1}=(1,0,0,\ldots,0)\in \R^d$ and let $\ol{B}=B_{\ol{1},1} \bigcup
B_{Y,\|Y\|}$. Introduce the random variable
\begin{equation}
\label{W}
   W= \frac{\lambda(\ol{B})}{\lambda(B_{0,1})}
\end{equation}
and let 
\[
  \alpha(d) \defeq \EXP\left[ \frac{2}{W^2}\right]~.
\]
The following result is proved in Section \ref{proofofthm1}.

\begin{theorem}
\label{thm:main}
Assume  that $\mu$ has a density $f$. Then
\begin{itemize}
\item[(i)]
\[
n \EXP\left[\mu(S_1)\mid X_1=x\right]
\to 1 \quad \text{for $\mu$-almost all $x$}~.
\]
\item[(ii)]
\[
n^2
\EXP\left[\mu(S_1)^2\mid X_1=x\right]
\to \alpha (d) \quad \text{for $\mu$-almost all $x$}~.
\]
\end{itemize}
\end{theorem}

In Section \ref{sec:alpha} we obtain estimates for the asymptotic
conditional second moment $\alpha(d)$. In particular, in Theorem \ref{thm:Luc4}
we show that for all dimensions, $1\le \alpha(d)\le 1+ 6(3/4)^{d/2}$
and therefore the asymptotic variance of $\mu(S_1)$ (conditioned on $X_1=x$)
decreases to zero exponentially in $d$. 
In the next result (Theorem \ref{thm:Luc3}) we determine the asymptotic
distribution of $\mu(S_1)$ (still conditioned on $X_1=x$).
We do this by determining the asymptotic moments of the limiting
distribution. Once again, the limit is the same for all $x$.

In order to describe the asymptotic moments, for any positive integer
define the random variable
\[
   W_k= \frac{\lambda(B_{\ol{1},1} \bigcup B_{Y_1,\|Y_1\|}\bigcup \ldots \bigcup B_{Y_{k-1},\|Y_{k-1}\|})}{\lambda(B_{0,1})}~,
\]
where $Y_1,\ldots , Y_{k-1}$ are independent 
random variables distributed uniformly in $B_{0,1}$.
Note that
\[
1=W_1\le W_2\le \cdots \le \frac{\lambda(B_{0,2})}{\lambda(B_{0,1})}=2^d.
\]
Now we may define a non-negative random variable $Z$ with moments
\[
\EXP[Z^k]=\EXP\left[ \frac{k!}{W_k^k}\right]
\]
for $k\ge 1$. We may use Carleman's condition to verify that the
distribution of $Z$ is uniquely defined. Indeed,
note that
\[
\EXP[Z^k]\le k!
\]
and therefore
\[
 \sum_{k=1}^{\infty}(\EXP[Z^{k}])^{-1/(2k)}\ge \sum_{k=1}^{\infty} (k!)^{-1/(2k)}=\infty~,
\]
and Carleman's condition is satisfied.
Note that if $E$ is an exponential $(1)$ random variable, then
\[
\EXP\left[E^k\right] = k! \ge \EXP[Z^k]\ge k!/2^{dk}=\EXP\left[\left(\frac{E}{2^d}\right)^k\right]~.
\]
We also have
\[
\EXP\left[e^{sZ}\right] \le \sum_{k=0}^{\infty}s^k=\frac{1}{1-s}
\]
for $0<s<1$ and
\[
\EXP\left[ e^{sZ}\right]\ge \sum_{k=0}^{\infty}\left(\frac{s}{2^d}\right)^k=\frac{1}{1-s/2^d}
\]
for $0<s<2^d$.

The next theorem establishes the convergence announced above.
The proof is sketched in Section \ref{proofofthm2}.

\begin{theorem}
\label{thm:Luc3}
Assume  that $\mu$ has a density $f$. Then, for $\mu$-almost all $x$, we have
that, conditionally on the event 
$X_1=x$, the random variable $n\mu(S_1)$ converges, in distribution,
to $Z$.
\end{theorem}

Note that for the case of a Voronoi tessallation of $\R^d$ defined by a Poisson point process
of constant intensity, Zuyev \cite{Zuy92} describes the distribution
of the volume of the so-called ``fundamental region'' of the cell containing the origin, conditionally on
having a point at the origin, as a mixture of Gamma
distributions. The fundamental region contains the Voronoi cell.
Since this distribution equals the limit for the
uniform density and our result is density free, the random variable $Z$ described here is stochastically dominated by the same
mixture of Gamma random variables.

In the case of $d=1$ it is easily seen that $Z$ is distributed as
$(E_1+E_2)/2$, where $E_1$ and $E_2$ are independent 
exponential$(1)$ random variables.

\section{Some values of $\alpha(d)$}
\label{sec:alpha}

In this section we investigate the asymptotic second moment $\alpha(d)$.
Since the limiting first moment equals $1$, we must have that 
$\alpha(d)\ge 1$. On the other hand, for all $d$, we have $\alpha(d)\le 2$.
To see this, recall from the proof of Theorem \ref{thm:main} that
\[
\alpha (d)/2 = \lim_{z\downarrow 0}\frac{\PROB\left\{\mu(B_{X,\|X-x \| }\cup B_{X',\|X'-x\| })\le z\right\}}{z^2}
\]
where $X,X'$ are i.i.d.\ with distribution $\mu$.
But clearly
\begin{align*}
&\PROB\left\{\mu(B_{X,\|X-x \| }\cup B_{X',\|X'-x \| })\le z\right\}\\
&=\PROB\left\{\mu(B_{X,\|X-x \| }\cup B_{X',\|X'-x \| })\le z,\mu(B_{X,\|X-x \| })\le z,\mu(B_{X',\|X'-x \| })\le z\right\}\\
&\le\PROB\left\{\mu(B_{X,\|X-x \| })\le z,\mu(B_{X',\|X'-x \| })\le z\right\} \\
&=\PROB\left\{\mu(B_{X,\|X-x \| })\le z\right\}^2
\end{align*}
and
\[
\frac{\PROB\left\{\mu(B_{X,\|X-x \| })\le z\right\}^2}{z^2}
\to 1~,
\]
once again from the proof of Theorem \ref{thm:main}.

It is not difficult to see that $\alpha(1) = 3/2$.
Indeed, by the definition (\ref{W}) of $W$,
\[
\ol{B}=B_{\ol{1},1} \cup
B_{Y,\|Y\|}=
B_{1,1} \cup
B_{Y,|Y|}
= \left\{  \begin{array}{ll}
           B_{1,1} & \mbox{if} \quad Y\ge 0 \\
          B_{1,1} \cup B_{Y,|Y|} & \mbox{if} \quad Y<0~,
          \end{array} \right.
\]
we have
\[
W= \frac{\lambda(\ol{B})}{\lambda(B_{0,1})}
\stackrel{\cal L}{=}
\left\{  \begin{array}{ll}
           1 & \mbox{ with probability } 1/2,\\
          1+U & \mbox{ with probability } 1/2,
          \end{array} \right.
\]
where $U$ is uniform $[0,1]$. Hence,
\[
\alpha(1)
= \EXP\left[ \frac{2}{W^2}\right]
= \frac 12 \left( \frac 21 + \EXP\left[ \frac{2}{(1+U)^2}\right] \right)
= 1 + \EXP\left[ \frac{1}{(1+U)^2}\right]
=3/2~.
\]

Previous work has considered the variance of the Lebesgue
measure of the Voronoi cell containing the origin defined by
a homogeneous Poisson process, conditioned on the fact that
a point falls in the origin. From these results, we deduce
values of $\alpha(d)$ for $d=2,3$. Indeed
Gilbert \cite{Gil62}, Brakke \cite{Bra85a}, and
Hayen and Quine \cite{HaQu02} showed that
$\alpha(2) \approx  1.2801760409267$ 
while Gilbert \cite{Gil62} and Brakke  \cite{Bra85b} showed
$\alpha(3)\approx  1.179032437845$.

Here we show that for large values of $d$, $\alpha(d)$ approaches $1$
exponentially fast. The proof is given in Section \ref{proofofthm3}.

\begin{theorem}
\label{thm:Luc4}
For all $d$,
\[
1\le \alpha(d)\le 1+ 6 (3/4)^{d/2}~.
\]
\end{theorem}

\section{On the diameter of a Voronoi cell}


Here we prove that, independently of the density, for $\mu$-almost all
$x\in \R^d$, conditionally on $X_1=x$, the diameter of the Voronoi cell
centered at $X_1$ converges to zero, in probability, at a rate of
$n^{-1/d}$.
More precisely, we have the following:

\begin{theorem}
\label{thm:diam}
Let $\mu$ have a density $f$. Then for $\mu$-almost all $x\in \R^d$,
\[
\lim_{t\to \infty}\limsup_{n\to \infty}\PROB\left\{n^{1/d}
  \diam (S_1)\ge t |X_1=x\right\}=0~.
\]
\end{theorem}

In particular, for $\mu$-almost all $x\in \R^d$, conditionally on $X_1=x$, $\diam(S_1)\to 0$
in probability. The theorem is proved in Section \ref{sec:diam}.

\section{Proofs}

The proofs use a version of the Lebesgue density theorem that we recall first.

We say that a class $\B$ of Borel sets in $\R^d$ is  \emph{good} if the following two conditions hold:
\[
\sup_{B\in \B}\frac{\lambda(\mbox{smallest ball }B_{0,r}\mbox{ containing }B)}{\lambda(B) }<\infty
\]
and
\[
\mbox{there exists a sequence }B_k\in \B \mbox{ with } \lambda(B_k)\downarrow 0.
\]
We say that $x$ is a \emph{Lebesgue point} for $f$ if for all good classes of Borel sets $\B$,
and all sequences $B_k\in \B$  with  $\lambda(B_k)\downarrow 0$,
\[
\lim_{k\to \infty}\frac{\int_{x+B_k}f}{\lambda(B_k)}=f(x)~.
\]
Let $A$ be the set of all $x\in\R^d$ such that $f(x)>0$ and $x$ is a Lebesgue point for $f$. Then $\mu(A)=1$
by Wheeden and Zygmund \cite[pp.\ 106--108]{WhZy77}. See also 
Devroye and Gy\"orfi \cite{DeGy85}, Chapter 2.

\subsection{Proof of Theorem \ref{thm:main}}
\label{proofofthm1}

\subsubsection*{Proof of part (i).}
Observe that
\begin{align*}
\EXP\left[\mu(S_1)\mid X_1=x\right]
&=
\PROB\left\{X_{n+1}\in S_1\mid X_1=x \right\}\\
&=
\PROB\left\{\cap_{i=2}^n\{X_i\notin B_{X_{n+1},\|X_{n+1}-x \| }\}\right\}\\
&=
\EXP\left[(1-Z(x))^{n-1} \right]~,
\end{align*}
where
\[
Z(x)=\mu(B_{X,\|X-x \| })~.
\]
It follows by integration by parts and the dominated convergence theorem
that
\[
n\EXP\left\{(1-Z(x))^{n-1} \right\}\to 1
\]
whenever 
\begin{equation}
\label{g1}
\lim_{z\downarrow 0}\frac{\PROB\left\{\mu(B_{X,\|X-x\| })\le z\right\}}{z}=1~.
\end{equation}
The intuitive reason of why such convergence should hold is that for
any $x$,
$\mu(B_{x,\|X-x\|})$ is uniformly distributed on $[0,1]$ and that
$\mu(B_{X,\|X-x\|}) \approx \mu(B_{x,\|X-x\|})$ when $\|X-x\|$ is
small. The rest of the proof establishes this convergence.

By the Lebesgue density theorem, it suffices to prove (\ref{g1}) for all 
Lebesgue points $x$ with $f(x)>0$. Fix such a point $x$.
Let $\B$ be the class of all closed balls of $\R^d$ containing the origin.
Since for any sequence $B_k\in \B$  with $\lambda(B_k)\downarrow 0$ we have
\[
\frac{\mu(x+B_k)}{\lambda(x+B_k)}\to f(x)~,
\]
for any $\epsilon \in (0,1)$ we can find  $\delta>0$ (possibly depending on  $x$)
such that $\|v-x\|\le \delta$ implies
\[
\left|\frac{\mu(B_{v,\|v-x\|})}{\lambda(B_{v,\|v-x\|})}- f(x)\right|\le \epsilon f(x)
\]
and
\[
\left|\frac{\mu(B_{x,\|v-x\|})}{\lambda(B_{x,\|v-x\|})}- f(x)\right|\le \epsilon f(x)~.
\]
This also implies that for any $v$ with $\|v-x\|\ge \delta$,
\begin{equation}
\label{g3}
\mu(B_{v,\|v-x\|})\ge \mu(B_{v^*,\delta})\ge (1-\epsilon)f(x)\lambda(B_{0,\delta})~,
\end{equation}
where $v^*$ is the unique point on the surface of $B_{x,\delta}$ and on the line segment $(x,v)$.
Take $z>0$ so small that
\[
z<(1-\epsilon)f(x)\lambda(B_{0,\delta}).
\]
Note also that
\[
\mu(B_{x,\|X-x\|})\stackrel{\cal L}{=}U~, 
\]
where $U$ is a uniform random variable on $[0,1]$.
We rewrite
\[
\mu(B_{X,\|X-x\|})
=\frac{\mu(B_{X,\|X-x\|})}{\lambda(B_{X,\|X-x\|})}\cdot \frac{\lambda(B_{x,\|X-x\|})}{\mu(B_{x,\|X-x\|})}\cdot  \mu(B_{x,\|X-x\|})~.
\]
If $\|X-x\|\le \delta$, then the first two factors are sandwiched between
\[
\frac{f(x)(1-\epsilon)}{f(x)(1+\epsilon)} \quad \mbox{ and } \quad \frac{f(x)(1+\epsilon)}{f(x)(1-\epsilon)}~.
\]
Since $\mu(B_{X,\|X-x\|})\le z$ implies $\|X-x\|\le \delta$ (see (\ref{g3})), we have
\begin{eqnarray*}
\PROB\{\mu(B_{X,\|X-x\|})\le z\}
&=&\PROB\{\mu(B_{X,\|X-x\|})\le z, \|X-x\|\le \delta\}\\
&\le&\PROB\left\{\frac{f(x)(1-\epsilon)}{f(x)(1+\epsilon)}\mu(B_{x,\|X-x\|})\le z, \|X-x\|\le \delta\right\}\\
&\le&\PROB\left\{\frac{1-\epsilon}{1+\epsilon}U\le z\right\}\\
&=&\min\left\{z\frac{1+\epsilon}{1-\epsilon},1\right\}~.
\end{eqnarray*}
Similarly,
\begin{align*}
\PROB\{\mu(B_{X,\|X-x\|})\le z, \|X-x\|\le \delta\}
&\ge\PROB\left\{\frac{1+\epsilon}{1-\epsilon}\mu(B_{x,\|X-x\|})\le z, \|X-x\|\le \delta,\mu(B_{X,\|X-x\|})\le z\right\}\\
&\ge\PROB\left\{\frac{1+\epsilon}{1-\epsilon}U\le z\right\}
-\PROB\left\{ \|X-x\|> \delta,\mu(B_{X,\|X-x\|})\le z\right\}\\
&=z\frac{1-\epsilon}{1+\epsilon}~.
\end{align*}
This proves (\ref{g1}) and part (i) of Theorem \ref{thm:main}.
\qed

\subsubsection*{Proof of part (ii).}
Similarly to the proof of part (i), observe that
\begin{align*}
\EXP\left[\mu(S_1)^2\mid X_1=x\right]
&=
\PROB\left\{X_{n+1}\in  S_1,X_{n+2}\in  S_1\mid X_1=x\right\} \\
&=
\PROB\left\{\cap_{i=2}^n\{X_i\notin B_{X_{n+1},\|X_{n+1}-x \| }\cup B_{X_{n+2},\|X_{n+2}-x \| }\}\right\}\\
&=
\EXP\left[(1-Z_2(x))^{n-1} \right]~,
\end{align*}
where
\[
Z_2(x)=\mu(B_{X,\|X-x \| }\cup B_{X',\|X'-x \| })
\]
with $X$ and $X'$ independent and distributed as $\mu$.
In analogy with the argument of part (i), in order to prove that
\[
\lim_{n\to \infty} n^2\EXP\left\{(1-Z_2(x))^{n-1} \right\} = \alpha (d)~,
\]
it suffices to show that
\[
\lim_{z\downarrow 0}\frac{\PROB\left\{\mu(B_{X,\|X-x \| }\cup B_{X',\|X'-x\| })\le z\right\}}{z^2}=\alpha (d)/2.
\]
The rough idea of the proof is as follows. The approximate equalities
are made rigorous below. For small $z$,
\begin{eqnarray*}
\lefteqn{
\frac{\PROB\left\{\mu(B_{X,\|X-x\|}\cup  B_{X',\|X'-x\|})\le
    z\right\}}{z^2}  } \\
&=&
\frac{\PROB\left\{\frac{\mu(B_{X,\|X-x\|}\cup  B_{X',\|X'-x\|})}{\max\{ \mu(B_{X,\|X-x\|}), \mu(B_{X',\|X'-x\|})\} }\max\{ \mu(B_{X,\|X-x\|}), \mu(B_{X',\|X'-x\|})\}\le z\right\}}{z^2}\\
&\approx &
\frac{\PROB\left\{\frac{\mu(B_{X,\|X-x\|}\cup  B_{X',\|X'-x\|})}{\max\{ \mu(B_{X,\|X-x\|}), \mu(B_{X',\|X'-x\|})\} }\max\{ \mu(B_{0,\|X\|}), \mu(B_{0,\|X'\|})\}\le z\right\}}{z^2}\\
&\approx &\frac{\PROB\left\{W \max\{ U_1, U_2\}\le z\right\}}{z^2}\\
& & \text{(where $U_1,U_2$ are i.i.d.\ uniform, independent of $W$)} \\
&= &\frac{\PROB\left\{W U^{1/2}\le z\right\}}{z^2}\\
& = & \frac{\EXP \left[ \min (z^2/W^2 , 1 ) \right]}{z^2} \approx \EXP\left[\frac{1}{W^2}\right]~.
\end{eqnarray*}
To prove the desired limit formally, as before, by the Lebesgue density theorem, we may assume that
$x\in\R^d$ is such that $f(x)>0$ and $x$ is a Lebesgue point for $f$.
A key point of the proof uses coupling. 
Let $(Y_1,Y_2)$ be the canonical reordering of $(X,X')$ such that
\[
\|Y_2-x\|\ge \|Y_1-x\|~.
\]
and introduce $M=\max(\|X-x\|,\|X'-x\|) = \|Y_2-x\|$.
Define the random variable $N$ by
\[
   N = \left\{ \begin{array}{ll}
               1 & \text{if $Y_1=X'$} \\
               2 & \text{if $Y_2=X'$}~. 
                   \end{array} \right.
\]
Then set $V_2=Y_2$ and let $V_1$ be uniformly distributed on $B_{x,\|V_2-x\|}$
such that $V_1$ is maximally coupled with $Y_1$ 
given $Y_2$.
From Doeblin's coupling argument,
\[
\PROB\{Y_1\neq V_1\mid Y_2 \}=\frac 12 \int |f_{Y_1}(v)-f_{V_1}(v)|dv~,
\]
where $f_{Y_1}$, $f_{V_1}$ are the conditional densities of $Y_1$ and $V_1$ given $Y_2$.

Choose $\delta>0$ so small that for $M\le \delta$, we have, simultaneously,
\[
\frac{\mu(B_{x,M})}{\lambda(B_{x,M})}\in [f(x)(1-\epsilon),f(x)(1+\epsilon)]~,
\]
\[
\frac{\mu(B_{X,\|X-x \| }\cup B_{X',\|X'-x \| })}{\lambda(B_{X,\|X-x \| }\cup B_{X',\|X'-x \| })}\in [f(x)(1-\epsilon),f(x)(1+\epsilon)]~,
\]
and
\[
\frac{\mu(B_{X,M})}{\lambda(B_{X,M})}\left|\frac{\lambda(B_{x,M})}{\mu(B_{x,M})}-\frac{1}{f(x)} \right|
+
\frac{1}{\lambda(B_{x,M})f(x)}\int_{B_{x,M} } |f(v)-f(x)|dv
\le \epsilon~.
\]
Such a $\delta$ exists by three applications of the Lebesgue density theorem. (Recall that $x$ is a Lebesgue point.)
Since
\[
f_{Y_1}(v)=\frac{f(v)}{\mu(B_{x,\|Y_2-x \|})}\IND_{v\in B_{x,M}}
\quad \text{and} \quad
f_{V_1}(v)=\frac{1}{\lambda(B_{x,M})}\IND_{v\in B_{x,M}}~,
\]
(where $\IND$ denotes the indicator function) we have, writing $B=B_{x,M}$,
\begin{eqnarray*}
\int |f_{Y_1}(v)-f_{V_1}(v)|dv
&=& \int_B \left|\frac{f(v)}{\lambda(B)}\frac{\lambda(B)}{\mu(B)} -\frac{1}{\lambda(B)}\right|dv\\
&\le& \frac{1}{\lambda(B)}\int_B f(v)\left|\frac{\lambda(B)}{\mu(B)} -\frac{1}{f(x)}\right|dv
+\frac{1}{\lambda(B)}\int_B \left|\frac{f(v)}{f(x)} -1\right|dv\\
&=& \frac{\mu(B)}{\lambda(B)} \left|\frac{\lambda(B)}{\mu(B)} -\frac{1}{f(x)}\right|
+\frac{1}{f(x)}\frac{1}{\lambda(B)}\int_B \left|f(v)-f(x)\right|dv\\
&\le& \epsilon
\end{eqnarray*}
if $M\le\delta$, by choice of $\delta$.
Finally, define a pair of random variables $(V,V')$, both taking values in
$\R^d$, as follows.
\[
(V,V')= \left\{  \begin{array}{ll}
          (V_1,V_2) & \mbox{if} \quad N=2\\
          (V_2,V_1) & \mbox{if} \quad N=1~,
          \end{array} \right.
\]
so that
\[
\PROB\{(V,V')\neq (X,X')\mid M\}\le \IND_{M>\delta}+\IND_{M\le\delta}\frac{\epsilon}{2}~.
\]
Since
\[
(\mu(B_{x,\|X-x\|}),\mu(B_{x,\|X'-x\|}))\stackrel{\cal L}{=}(U,U'),
\]
where $U,U'$ are independent uniform  $[0,1]$ random variables, we have,
\[
\mu(B_{x,M})\stackrel{\cal L}{=}\max(U,U')\stackrel{\cal L}{=}\sqrt{U}~.
\]
By construction, $V_1$ is uniform on $B_{x,\|Y_2-x\|}$, so that, given $Y_2$,
\[
\frac{\lambda(B_{Y_2,\|Y_2-x \| }\cup B_{Y_1,\|Y_1-x \| })}{\lambda(B_{Y_2,\|Y_2-x \| })}
\stackrel{\cal L}{=}W~,
\]
where $W$ was defined in (\ref{W}).
To complete the argument, set
\[
B_X=B_{X,\|X-x \| },\quad B_{X'}=B_{X',\|X'-x \| },\quad M=\max(\|X-x \|,\|X'-x \|)~.
\]
Then
\begin{align*}
\mu(B_X\cup B_{X'})
&=\frac{\mu(B_X\cup B_{X'}) }{\lambda(B_X\cup B_{X'}) }\cdot \frac{\lambda(B_X\cup B_{X'}) }{\lambda(B_{x,M}) }
\cdot \frac{\lambda(B_{x,M}) }{\mu(B_{x,M}) }\cdot \mu(B_{x,M})\\
&\defeq I\cdot II\cdot III\cdot IV.
\end{align*}
Note that
\[
I\in [f(x)(1-\epsilon),f(x)(1+\epsilon)]
\]
when $M\le \delta$, and similarly,
\[
III\in \left[\frac{1}{f(x)(1+\epsilon)},\frac{1}{f(x)(1-\epsilon)}\right]
\]
when $M\le \delta$. When $(X,X')=(V,V')$, we have
\[
II=\frac{\lambda(B_X\cup B_{X'}) }{\lambda(B_{x,M}) }\stackrel{\cal L}{=}W
\]
with $W$ independent of
\[
IV=\mu(B_{x,M})\stackrel{\cal L}{=}\sqrt{U}~.
\]
Thus, since for small enough $z$, $\mu(B_X\cup B_{X'})\le z $ implies $M\le\delta$ (as argued in the proof of (\ref{g1})), for such $z$, we have
\[
\PROB\{\mu(B_X\cup B_{X'})\le z,(X,X')=(V,V') \}\le \PROB\left\{\frac{1-\epsilon}{1+\epsilon}W\sqrt{U}\le z \right\}
\]
and thus,
\[
\PROB\{\mu(B_X\cup B_{X'})\le z \}
= \PROB\{\mu(B_X\cup B_{X'})\le z, M\le \delta \}+\PROB\{\mu(B_X\cup B_{X'})\le z, M>\delta\}~.
\]
Clearly,
\[
\PROB\{\mu(B_X\cup B_{X'})\le z, M>\delta\}=0
\]
for $z$ small enough. For such a $z$, we have
\begin{eqnarray*}
 \PROB\{\mu(B_X\cup B_{X'})\le z, M\le \delta \}
&=& \PROB\{\mu(B_X\cup B_{X'})\le z, M\le \delta,(X,X')\neq (V,V') \}\\
& & \quad +\PROB\{\mu(B_X\cup B_{X'})\le z, M\le \delta,(X,X')=(V,V') \}\\
&\defeq & I+ II.
\end{eqnarray*}
We have
\begin{align*}
I
&\le
\PROB\{\mu(B_X)\le z\}
\PROB\{\mu(B_{X'})\le z\}
\sup_{\rho\le\delta}\PROB\{(X,X')\neq (V,V')\mid M\le \rho \}\\
&=
z^2(1+o(1))
\sup_{\rho\le\delta}\PROB\{(X,X')\neq (V,V')\mid M\le \rho \}\\
&\qquad \qquad \mbox{(by the proof of (\ref{g1}))}\\
&=
z^2(1+o(1))\epsilon \qquad \mbox{(by the choice of  $\delta$)}~.
\end{align*}
Also,
\begin{align*}
II
&\le
\PROB\{\mu(B_X\cup B_{X'})\le z, M\le \delta \}\\
&\le
\PROB\left\{\frac{1-\epsilon}{1+\epsilon}W\sqrt{U}\le z \right\}\\
&\le
\left(\frac{1+\epsilon}{1-\epsilon}\right)^2\EXP\left\{\frac{1}{W^2} \right\}z^2\\
&=
\left(\frac{1+\epsilon}{1-\epsilon}\right)^2z^2\alpha(d)/2~.
\end{align*}
On the other hand,
\begin{align*}
&\PROB\{\mu(B_X\cup B_{X'})\le z \}\\
&\ge\PROB\{\mu(B_X\cup B_{X'})\le z, M\le \delta, (X,X')= (V,V') \}\\
&\ge\PROB\left\{\frac{1+\epsilon}{1-\epsilon}\frac{\lambda(B_V\cup B_{V'})}{\lambda(B_{x,M})}\mu(B_{x,M})\le z,\mu(B_X\cup B_{X'})\le z, M\le \delta, (X,X')= (V,V') \right\}\\
&\ge
\PROB\left\{\frac{1+\epsilon}{1-\epsilon}W\sqrt{U}\le z \right\}
 -\PROB\left\{\mu(B_X\cup B_{X'})\le z, M\le \delta, (X,X')\neq (V,V') \right\}\\
&\ge
\PROB\left\{U\le \left(\frac{z(1-\epsilon)}{W(1+\epsilon)}\right)^2 \right\}
 -\PROB\{\mu(B_X)\le z,\mu(B_{X'})\le z,M\le \delta, (X,X')\neq (V,V') \}~,
\end{align*}
and therefore
\begin{align*}
&\PROB\{\mu(B_X\cup B_{X'})\le z \}\\
&\ge
\EXP\left\{\min\left\{\left(\frac{z(1-\epsilon)}{W(1+\epsilon)}\right)^2,1 \right\}\right\}\\
&\quad
-\PROB\{\mu(B_X)\le z\}
\PROB\{\mu(B_{X'})\le z\}
\sup_{\rho\le\delta}\PROB\{(X,X')\neq (V,V')\mid M\le \rho \}\\
&=
(1+o(1))z^2\left(\frac{1-\epsilon}{1+\epsilon}\right)^2
\EXP\left\{\frac{1}{W^2}\right\}\\
&\quad
-(1+o(1))z^2
\sup_{\rho\le\delta}\PROB\{(X,X')\neq (V,V')\mid M\le \rho \}\\
&\qquad \qquad \mbox{(by the dominated convergence theorem, and proof of (\ref{g1}))}\\
&\ge
(1+o(1))\left(\frac{1-\epsilon}{1+\epsilon}\right)^2z^2\alpha(d)/2-(1+o(1))z^2\epsilon\\
&\qquad \qquad \mbox{(by the choice of $\delta$)}.
\end{align*}
Since $\epsilon$ was arbitrary, we are done.
\qed

\subsection{Sketch of proof of Theorem \ref{thm:Luc3}}
\label{proofofthm2}

Since the proof of Theorem \ref{thm:Luc3} 
is an extension of that of Theorem \ref{thm:main}, we only sketch 
the arguments.

By the moment method, it suffices to show that for 
all Lebesgue points $x\in\R^d$ with $f(x)>0$, and for all $k\ge 1$,
we have
\[
\EXP\left[(n\mu(S_1))^k|X_1=x\right]\to \EXP[Z^k]~,
\]
As we argued in the case $k=2$ in the proof of Theorem \ref{thm:main},
\[
\EXP\left[n^k\mu(S_1))^k|X_1=x\right]
=
n^k\EXP\left[(1-Z_k(x))^{n-1} \right],
\]
where
\[
Z_k(x)\defeq \mu(B_{X_1,\|X_1-x\|}\cup \cdots \cup B_{X_{k},\|X_{k}-x\|})
\stackrel{\cal L}{\approx}
W_kU^{1/k}
\]
where $U$ is uniform $[0,1]$. Here we use the fact that
\[
\max_{1\le i\le k} \mu(B_{0,\|X_i\|}) \stackrel{\cal L}{=}\max_{1\le i\le k}U_i\stackrel{\cal L}{=}U^{1/k}
\]
with the $U_i$ being independent and uniform on $[0,1]$.

Just like in the proof of Theorem \ref{thm:main}, in order to show
that
\[
\EXP\left[n^k\mu(S_1))^k|X_1=x\right] \to \EXP\left[ \frac{k!}{W_k^k}\right]~,
\]
it suffices to show that
\[
\lim_{z\downarrow 0}\frac{\PROB\left\{\mu(B_{X_1,\|X_1-x\|}\cup \cdots \cup B_{X_{k},\|X_{k}-x\|})\le z\right\}}{z^k}=\EXP\left[\frac{1}{W_k^k}\right]~.
\]
By the approximation above, for small $z$,
\[
\frac{\PROB\left\{\mu(B_{X_1,\|X_1-x\|}\cup \cdots \cup B_{X_{k},\|X_{k}-x\|})\le z\right\}}{z^k} \approx 
\frac{\PROB\left\{W_kU^{1/k}\le z\right\}}{z^k} \approx \EXP\left[\frac{1}{W_k^k}\right]~.
\]
The approximation can be made rigorous by the same arguments as detailed in
the proof of Theorem \ref{thm:main}.
\qed

\subsection{Proof of Theorem \ref{thm:Luc4}}
\label{proofofthm3}

Define
\[
A=B_{0,1},\quad B=B_{\ol{1},1}, \quad \mbox{ and } \quad C=B_{Y,\|Y\|}~,
\]
where $Y$ is uniformly distributed on $A$.
Define
\[
W=\frac{\lambda(B\cup C)}{\lambda(B) } \quad \mbox{ and }\quad U=\frac{\lambda(C)}{\lambda(B) }
\]
and observe that $U$ is uniformly distributed on $[0,1]$.
We write
\begin{align*}
\alpha(d)-1
&= \EXP\left[ \frac{2}{W^2}-   \frac{2}{(1+U)^2}\right]\\
&= 2\EXP\left[ \frac{1}{W^2}\left(1-   \left(\frac{W }{1+U}\right)^2\right)\right] \\
&\le  2\EXP\left[ 1-   \left(\frac{W }{1+U}\right)^2\right]~,
\end{align*}
since $W\ge 1$. We have that
\[
2\left[ 1-   \left(\frac{W }{1+U}\right)^2\right]
\le 2
\]
and
\begin{align*}
2\left[ 1-   \left(\frac{W }{1+U}\right)^2\right] 
&=  2\left[ \frac{(1+U-W)(1+U+W) }{(1+U)^2}\right]\\ 
&\le  4\left[ \frac{1+U-W }{1+U}\right]\\
&= 4\left[ \frac{\lambda(B\cap C)}{\lambda(B)+\lambda(C)}\right]~.
\end{align*}
Thus
\begin{align*}
\alpha(d)-1
&\le  
2\EXP\left[\IND_{Y\in B}\right]
+ 4\EXP\left[ \frac{\lambda(B\cap C)}{\lambda(B)+\lambda(C)}\IND_{Y\notin B}\right]~.
\end{align*}
We finish the proof by showing that
\begin{equation}
\label{c1}
\EXP\left[ \IND_{Y\in B}\right] \le (3/4)^{d/2} 
\end{equation}
and
\begin{equation}
\label{c2}
\EXP\left[ \frac{\lambda(B\cap C)}{\lambda(B)+\lambda(C)}\IND_{Y\notin B}\right] \le (3/4)^{d/2}~.
\end{equation}
For (\ref{c1}), note that
\[
\EXP\left[\IND_{Y\in B}\right]=\frac{\lambda(A\cap B)}{\lambda(B)}\le \frac{\lambda( B_{b,\sqrt{3/4}})}{\lambda(B)}=(3/4)^{d/2}~,
\]
where $b=(1/2,0,0,\ldots 0)\in \R^d$,
since $A\cap B\subset B_{b,\sqrt{3/4}}$.
For (\ref{c2}), we bound
\begin{align*}
\frac{\lambda(B\cap C)}{\lambda(B)+\lambda(C)}\IND_{Y\notin B}
&\le\sup_{y\notin B}
\frac{\lambda(B\cap B_{y,\|y\|})}{\lambda(B)+\lambda(B_{y,\|y\|})}\\
&\le
\frac{\lambda(B \cap B_{a,1})}{\lambda(B)}\\
& \quad \text{(where $a=(1/2,\sqrt{3/4},0,0,\ldots ,0)\in \R^d$)} \\
&=
\frac{\lambda(A\cap B)}{\lambda(B)}\\
&\le
(3/4)^{d/2}
\end{align*}
by arguing as above. That the supremum reached by placing $Y$ (thus, $y$) at $a$ is clear in two steps.
First, the intersection can only grow by replacing $y$ by $y/\|y\|$ since
\[
B_{y,\|y\|}\subset B_{y/\|y\|,1}.
\]
Next, of all the points on the surface, but outside $B_{\ol{1},1} $, the intersection $\lambda(B\cap B_{y/\|y\|,1}) $
is maximized by placing $y$ at $a$.
\qed

\subsection{Proof of Theorem \ref{thm:diam}}
\label{sec:diam}


Let $\gamma_d$ be the minimal number of cones $C_1,\ldots,C_{\gamma_d}$ of angle $\pi/4$ 
centered at $0$ such that their union covers $\R^d$. 
Let $R_{n,j}$, $j=1,\ldots ,\gamma_d$, be the distance between $X_1$
and the nearest neighbor among $X_2,\ldots ,X_n$ belonging to
$X_1+C_j$ (i.e., the cone $C_j$ translated by $X_1$).
Define $R_{n,j}=\infty$ if no such point exists.

We bound the diameter of he Voronoi cell $S_1$ by observing that
\[
\diam(S_1)\le \sqrt{d}\max_{j=1,\ldots ,\gamma_d}R_{n,j}~.
\]
A simple extension of the Lebesgue density theorem implies that 
if $B=B_{0,1}$ is the unit ball centered at the origin, then for $\mu$-almost all $x\in \R^d$,
\begin{equation}
\label{*}
\min_{j=1,\ldots ,\gamma_d}\frac{\int_{x+r[C_j\cap
    B]}f}{\lambda(r[C_j\cap B])}\to f(x) \quad \mbox{ as } \quad r\downarrow 0~.
\end{equation}
Thus, for $\mu$-almost all $x$, there exists $R(x)>0$ such that for all $0<r\le R(x)$,
\[
\min_{j=1,\ldots ,\gamma_d}\int_{x+r[C_j\cap B]}f \ge r^d \frac{f(x)}{2} \lambda(C_1\cap B)~.
\]
If $f(x)=0$ or $x$ does not satisfy (\ref{*}), set $R(x)=0$. For any
$t>0$, we have
\begin{align*}
\left\{\diam(S_1)>tn^{-1/d}\right\} 
& \subset
\left\{\max_{j=1,\ldots ,\gamma_d}R_{n,j}>\frac{tn^{-1/d}}{\sqrt{d}}\right\}\\
&\subset
\bigcup_{j=1}^{\gamma_d}\left\{X_1+(C_j\cap B)\frac{tn^{-1/d}}{\sqrt{d}} \mbox{ has no point among } X_2,\ldots ,X_n\right\}~.
\end{align*}
Thus, we have
\begin{eqnarray*}
\lefteqn{
\PROB\left\{n^{1/d}  \diam (S_1)\ge t |X_1=x\right\}   } \\
& \le & \sum_{j=1}^{\gamma_d}
\PROB\left\{x+(C_j\cap B)\frac{tn^{-1/d}}{\sqrt{d}} \mbox{ has no
    point among } X_2,\ldots ,X_n\right\}~.
\end{eqnarray*}
We bound the probability of each event in the
union as follows.
\begin{eqnarray*}
\lefteqn{
\PROB\left\{x+(C_j\cap B)\frac{tn^{-1/d}}{\sqrt{d}} \mbox{ has no
    point among } X_2,\ldots ,X_n\right\}  } \\
&\le &\PROB\left\{x+(C_j\cap B)\min\left(R(x),\frac{tn^{-1/d}}{\sqrt{d}}\right) \mbox{ has no point among } X_2,\ldots ,X_n\right\}\\
&= &
\left(1-  \mu\left( x+(C_j\cap B)\min\left(R(x),\frac{tn^{-1/d}}{\sqrt{d}}\right)\right)\right)^{n-1}\\
&\le &
\left(1-  \left( \min\left(R(x),\frac{tn^{-1/d}}{\sqrt{d}}\right)
  \right)^d  \lambda(C_1\cap B) \frac{f(x)}{2}\right)^{n-1}
\\
& \le &
\exp\left(- (n-1) \min\left(R(x)^d,\frac{t^dn^{-1}}{\sqrt{d^d}}\right)\left(\lambda(C_1\cap B) \frac{f(x)}{2}\right)
\right)
\end{eqnarray*}
and the theorem follows since $R(x)^df(x)>0$ for $\mu$-almost all $x$.
\qed

\end{document}